# Inverse iterative simulation: An efficient approach for contaminant source identification

Jiangjiang Zhang[1]

[1] College of Environmental and Resource Sciences, Zhejiang University, Hangzhou, 310058, China

**Abstract.** In groundwater contaminant remediation and risk assessment, it is important to identify parameters of the contaminant source and hydraulic conductivity field by solving an inverse problem. However, if the dimensionality of the inverse problem is high, it is usually computationally expensive to obtain accurate estimation and uncertainty assessment of these parameters. This is particularly the case when Markov Chain Monte Carlo (MCMC) sampling is used. In this paper, an efficient approach entitled inverse iterative simulation (iIS) is proposed to efficiently identify the contaminant source characteristics, together with the hydraulic conductivity field. The iIS algorithm utilizes a simple approach borrowed from Ensemble Smother (ES) to update model parameters and an inverse Gaussian process (iGP) approach to improve the accuracy of parameter updating. Two numerical experiments are tested. For the low dimensional case (with 11 parameters), the iIS algorithm can obtain parameter estimation very close to that of MCMC method. For the high dimensional case (with 108 parameters), the iIS algorithm can obtain accurate parameter estimation with very low computational cost.

# 1. Introduction

For better prediction of the effect of human activities on subsurface environment, it is vital to develop accurate groundwater models. However, uncertainties derived from measurements and model structures are ubiquities. Moreover, limited data that rarely contains sufficient information to identify the subsurface characteristics further undermine accurate modeling [*Tartakovsky*, 2013]. Thus, uncertainty quantification is essential in groundwater modeling, where quantifying parametric uncertainty is the basis of quantifying model structure uncertainty [*Zhang et al.*, 2013].

In groundwater contaminant remediation and risk assessment, identifying parameters of contaminant source (e.g., source location and release history) and hydraulic conductivity field is essential. However, directly measuring these parameters is difficult or even impossible. Thus, we need to estimate the model parameters indirectly from concentration and hydraulic head measurements by solving an inverse problem. Many inverse methods have been used to identify containment source parameters, e.g., geostatistical approach [*Snodgrass and Kitanidis*, 1997; *Sun*, 2007], minimum relative entropy method [*Woodbury et al.*, 1998], correlation coefficient optimization [*Sidauruk et al.*, 1998], least squares methods [*Liu and Ball*, 1999], Genetic Algorithm [*Mahinthakumar and Sayeed*, 2005], simulated annealing [*Yeh et al.*, 2007], and Markov Chain Monte Carlo (MCMC) sampling [*Wang and Jin*, 2013; *Zeng et al.*, 2012; *Zhang et al.*, 2015]. Nowadays, MCMC methods are becoming increasingly popular in hydrologic model uncertainty quantification for its general applicability in highly nonlinear and non-Gaussian problems involving complex processes. However, for high dimensional problems, even with some advanced MCMC algorithms (e.g., DRAM [*Haario et al.*, 2006] and DREAM$_{(ZS)}$ [*Vrugt et al.*, 2008; *Vrugt et al.*, 2009]), a very large number of model evaluations are usually needed to sufficiently explore the posterior parameter space. When the uncertainty of contaminant source (source location and release history) and hydraulic conductivity field are considered simultaneously, the number of unknown parameters would be rather large, which would require very huge computational cost if MCMC method is adopted. To efficiently infer these parameters, in this paper, a new approach is proposed as an alternative to MCMC method. This method utilizes a simple approach borrowed from Ensemble Smother (ES) [*Evensen and Van Leeuwen*, 2000] to update model parameters and an inverse Gaussian process (iGP) approach to improve the accuracy of parameter updating. We have called this algorithm inverse iterative simulation (iIS). As a benchmark, the algorithm will be compared with the widely used MCMC algorithm DREAM$_{(ZS)}$ [*Vrugt et al.*, 2008; *Vrugt et al.*, 2009], which have shown great efficiency and accuracy in hydrologic model parameter inference. The paper continues with descriptions of contaminant transport model and the algorithm, followed by implementing the iIS algorithm on two synthetic examples, and ends with some conclusions.

## 2. Theory and Methods

### 2.1. Flow and Transport Model

In this study, the transport of nonreactive contaminant in a two-dimensional (2-D) heterogeneous flow field is considered. The steady state saturated groundwater flow satisfies the following governing equation [*Harbaugh*, 2005]:

$$\frac{\partial}{\partial x}\left(K_{xx}\frac{\partial h}{\partial x}\right)+\frac{\partial}{\partial y}\left(K_{yy}\frac{\partial h}{\partial y}\right)+W=S\frac{\partial h}{\partial t}, \qquad (1)$$

where $K_{xx}$ and $K_{yy}$ are values of hydraulic conductivity along the $x$ and $y$ coordinate directions $[LT^{-1}]$, with the assumption that they are parallel to the major axes of hydraulic conductivity; $h$ is the hydraulic pressure head $[L]$; $W$ is the sink or source term $[T^{-1}]$; $S$ is the specific storage of the porous material $[L^{-1}]$; $t$ is time $[T]$.

With appropriate initial and boundary conditions, the 2-D groundwater flow problem can be solved numerically. Then the transport of a nonreactive contaminant

can be obtained by solving the following advection dispersion equation [*Zheng and Wang*, 1999]:

$$\frac{\partial(\theta C)}{\partial t} = \frac{\partial}{\partial x_i}\left(\theta D_{ij} \frac{\partial C}{\partial x_j}\right) - \frac{\partial}{\partial x_i}(\theta v_i C) + q_s C_s, \qquad (2)$$

where $\theta$ is the porosity of the subsurface medium; $C$ is the dissolved concentration of contaminant [$ML^{-3}$]; $x_{i,j}$ are the distances along the respective Cartesian coordinate axes [L]; $D_{ij}$ is the hydrodynamics dispersion tensor [$L^2T^{-1}$]; $v_i$ is the seepage or linear pore water velocity [$LT^{-1}$]; $q_s$ is volumetric flow rate per unit volume of aquifer representing fluid sources (positive) [$T^{-1}$]; $C_s$ is the concentration of the source [$ML^{-3}$]. In this study, $q_s C_s$ [$ML^{-3}T^{-1}$] is treated as a single variable, which characterizes the contaminant source strength $S_s$. The above equations become stochastic when the conditions or parameters (or both) are uncertain.

## 2.2. Methods

In a hydrologic model, measurements **d** can be expressed as

$$\mathbf{d} = f(\mathbf{m}) + \boldsymbol{\varepsilon}, \qquad (3)$$

where **m** and $f(\mathbf{m})$ are $n_m \times 1$ and $n_d \times 1$ vectors of the model parameters and outputs, respectively, $n_m$ and $n_d$ are the dimensions of parameters and measurements, respectively; $\boldsymbol{\varepsilon}$ is a $n_d \times 1$ vector of measurement errors. We are interested in the estimation and uncertainty assessment of model parameters **m** from noisy measurements **d**. In this paper, we propose a simple while effective method which combines the updating scheme similar to that of Ensemble Smother and inverse Gaussian process to estimate the model parameters. The main processes are illustrated as follows.

(1) This method starts with drawing *N* samples from prior distribution of parameters

$\mathbf{M} = [\mathbf{m}_1, \mathbf{m}_2, ..., \mathbf{m}_N]$, then calculating their corresponding model outputs $\mathbf{F} = [f(\mathbf{m}_1), f(\mathbf{m}_2), ..., f(\mathbf{m}_N)]$, which can be easily realized in a parallel mode.

(2) With available measurements $\mathbf{d}$, the parameter samples can be updated with a scheme similar to that of Ensemble Smother, and the updated parameter samples are demoted as $\mathbf{M}^a$,

$$\mathbf{M}^a = \mathbf{M} + \mathbf{K}(\mathbf{En_d} - \mathbf{F}), \tag{4}$$

where $\mathbf{K}$ is the Kalman gain; $\mathbf{En_d}$ is the ensemble of perturbed measurements $[\mathbf{d}_1, \mathbf{d}_2, ..., \mathbf{d}_N]$, $\mathbf{d}_i = \mathbf{d} + \mathbf{\varepsilon}_i$, $\mathbf{\varepsilon}_i$ is one realization of measurement noise.

The Kalman gain $\mathbf{K}$ is calculated using the following equations,

$$\mathbf{K} = \mathbf{P_M}(\mathbf{P_F} + \mathbf{R})^{-1}, \tag{5}$$

$$\mathbf{P_M} = \frac{(\mathbf{M} - \overline{\mathbf{M}})(\mathbf{F} - \overline{\mathbf{F}})^T}{N-1}, \tag{6}$$

$$\mathbf{P_F} = \frac{(\mathbf{F} - \overline{\mathbf{F}})(\mathbf{F} - \overline{\mathbf{F}})^T}{N-1}, \tag{7}$$

where $\overline{\mathbf{M}}$ and $\overline{\mathbf{F}}$ are matrixes with $N$ columns, each column of the two matrixes are the mean of $\mathbf{M}$ and $\mathbf{F}$, respectably; $\mathbf{R}$ is the covariance matrix of measurement error.

(3) Set $\mathbf{M} = \mathbf{M}^a$, calculate the corresponding model outputs $\mathbf{F} = [f(\mathbf{m}_1), f(\mathbf{m}_2), ..., f(\mathbf{m}_N)]$.

(4) Repeat steps (2) and (3), until the stop criterion is satisfied. The stop criterion is defined according to the consistency between residuals (the difference between the latest $\mathbf{F}$ and $\mathbf{d}$) and measurement error statistics. This is realized as follows.

Given measurement $\mathbf{d}$ and measurement error distribution $N(\mathbf{0}, \mathbf{\sigma}^2)$, calculate the Gaussian likelihood values of $\mathbf{F}$, which are denoted as $Lik_i^{\mathbf{F}}, i = 1, 2, ..., N$. Meanwhile, the Gaussian likelihood values of $\mathbf{En_d}$ are also calculated, $Lik_i^{\mathbf{En_d}}, i = 1, 2, ..., N$.

Remove the outliers in $Lik_i^{\mathbf{F}}, i = 1, 2, ..., N$, if most (90%) of $Lik_i^{\mathbf{F}}$ are within the bounds of $Lik_i^{\mathbf{En_d}}$, then the iteration procedure stops. In practice, the log values of the likelihood are used.

**(5)** To improve the accuracy of the updating process described in step (2), the parameter samples in **M** are refined with an inverse Gaussian process ahead of updating. The brief idea is simple. Mapping from model outputs to parameters, Gaussian process regression is used to construct an inverse surrogate system of the original function. Given measurements **d** as inputs to this inverse surrogate system, parameter estimation $\mathbf{m}_{est}$ can be obtained directly. Then the sample **m'** in **M** and the corresponding $f(\mathbf{m'})$ in **F** with the smallest Gaussian likelihood value are replaced with $\mathbf{m}_{est}$ and $f(\mathbf{m}_{est})$, respectively. Totally, the 50 worst samples in **M** and **F** (i.e., with smaller likelihood values) are replaced step by step in this way.

## 3. Numerical Experiments

## 3.1. Case Study 1: Contaminant Source Identification with Zonated Conductivity Field

In this case study, we tested the iIS algorithm for a contaminant source identification problem in steady saturated flow.

As shown in Figure 1, the flow domain is 20[L] in *x* direction and 10[L] in *y* direction. In this case, the porosity was $\theta = 0.25$, the longitudinal dispersivity and the transverse dispersivity were $\alpha_L = 0.3 [L^2T^{-1}]$ and $\alpha_T = 0.03 [L^2T^{-1}]$, respectively. The conductivity field was represented with three hydraulic zones. In each zone, the hydraulic conductivity $K_i$ [LT$^{-1}$] (represented by its log value, i.e., $Y_i = \log K_i, i = 1, 2, 3$) was assumed to be homogenous with unknown values. With no flow boundaries in the upper and lower sides, constant head boundaries with pressure heads of 12[L] and 11[L] at the left and right sides, respectively, the flow equation was solved numerically with MODFLOW [*Harbaugh*, 2005]. In this steady saturated flow field, a contaminant source located within a potential area denoted by the red dashed rectangle in Figure 1 was released from 1[T] to 6[T]. Then the solute transport equations was solved numerically with MT3DMS [*Zheng and Wang*, 1999].

[Figure 1]

In this case, there were totally 11 unknown parameters, i.e., 3 conductivity parameters ($Y_1, Y_2, Y_3$), 8 contaminant source parameters, including the location ($x_s, y_s$) and time-varying strengths $S_{si}$ [MT$^{-1}$] for $t_i = i$ [T] : $(i+1)$ [T], $i = 1, 2, ..., 6$. Their distributions were assumed to be uniform with given ranges, as listed in Table 1. To solve this inverse problem, concentration and hydraulic head measurements at 5 locations (the blue dots in Figure 1) were used. Every 2[T] from $t = 2$[T] to $t = 10$[T], the concentration measurements were collected. Since the flow was steady, the hydraulic head measurements were sampled only once at the 5 locations. The errors for the concentration and head measurements were assumed to follow $N(0, 0.05^2)$ and $N(0, 0.01^2)$, respectively. With these noisy measurements, the inverse problem was solved with the iIS and DREAM$_{(ZS)}$ algorithms, respectively.

For the iIS algorithm, in each iteration, the number of parameter samples for updating was 400, and the 50 parameter samples with the smallest likelihood values were replaced with the inverse Gaussian process method as described in step (5), 2.2., which means that the total number of model evaluations was 450 in each iteration. Figure 2 shows the log Gaussian likelihood values of $\mathbf{En_d}$ (the ensemble of perturbed measurements) and $\mathbf{F}$ (the ensemble of model responses corresponding to the 400 parameter samples) at each iteration. At the first 5 iterations, the log likelihood values of $\mathbf{F}$ are far beyond the range of $\log Lik^{\mathbf{En_d}}$, which means that the distance between

**F** and measurements **d** is large. At the 6th iteration, most of the log likelihood values of **F** are within the bounds of $\log Lik^{En_d}$, which means that the residuals between **F** and **d** are consistent with the measurement error statistics. In other words, the number of model evaluations for the iIS algorithm is 6×450 = 2,700. The trace plots of the 11 parameters (the 50 parameter samples replaced in each iteration are not shown) are shown in Figure 3. It is obvious that, the parameter samples converge to the true parameter values step by step. As in each iteration, the calculation of model outputs **F** given parameter ensemble **M** can be easily realized in parallel, the time needed for the iIS algorithm is affordable.

[Figure 2]

[Figure 3]

To provide references of posterior parameter distributions, the DREAM$_{(ZS)}$ algorithm was implemented with three parallel chains and altogether 30,000 model evaluations, where Gaussian likelihood with heteroscedastic measurement errors was adopted. The convergence was reached after about 20,000 model evaluations, and the last 10,000 samples were used to estimate the posterior distributions of parameters. As shown in Figure 4, compared with DREAM$_{(ZS)}$, the iIS algorithm can obtain similar distributions of the 11 parameters. However, the distributions of some parameters (i.e., $S_{s2}$, $Y_2$ and $Y_3$) obtained by the iIS algorithm are with slightly higher peaks. This may be caused by the fact that, although we try to make sure that the residuals between **F** and measurements **d** are consistent with measurement error statistics, it is unavoidable to over update some parameters slightly.

[Figure 4]

The inverse Gaussian process described in step (5), 2.2. to refine parameter samples could guarantee the accuracy of the updating process described in step (2), if it is not applied, it may result in less accurate parameter estimation. To illustrate this, one another set of parameters randomly drawn from prior distributions were chosen as the true parameters, and the measurements were generated with additive measurement errors. Figure 5 shows the posterior distributions obtained by the DREAM$_{(ZS)}$ algorithm, the iIS algorithm adopting the inverse Gaussian process and the iIS algorithm without the inverse Gaussian process, respectively. It is clearly shown that, if the inverse Gaussian process is not used, the posterior distribution of parameters obtained by the iIS algorithm are more likely to deviate from those of MCMC.

[Figure 5]

## 3.2. Case Study 2: Contaminant Source Identification with Continuous Random Conductivity Field

In Case 1, the conductivity field with only three hydraulic zones is considered, which is an over-simplification of real cases. To be more realistic, the conductivity can be modeled as a continuously varying random field. In this case study, the log conductivity field $Y(\mathbf{x})$ is assumed as a spatially correlated Gaussian random field with separable exponential correlation form shown in Eq. (8)

$$C_Y(\mathbf{x}_1, \mathbf{x}_2) = C_Y(x_1, y_1; x_2, y_2) = \sigma_Y^2 \exp\left[-\frac{|x_1 - x_2|}{\lambda_x} - \frac{|y_1 - y_2|}{\lambda_y}\right], \quad (8)$$

where $\sigma_Y^2 = 1$ is the variance and $\lambda_x = 20[L]$ and $\lambda_y = 10[L]$ are correlation lengths in the *x* and *y* directions, respectively. In this case however, the number of unknown parameters would be very large, which poses unaffordable computational burden for parameter estimation. To alleviate this problem, some dimension reduction techniques can be employed to reduce the number of unknown parameters. In this case, the Karhooven Loeve (KL) expansion [*Zhang and Lu*, 2004] is used to reduce the number of unknown parameters from the total grid number 3,321 to 100 truncated KL terms,

$$Y(\mathbf{x}) \approx \bar{Y}(\mathbf{x}) + \sum_{i=1}^{100} \xi_i \sqrt{\lambda_i} f_i(\mathbf{x}), \quad (9)$$

where $\bar{Y}(\mathbf{x})$ is the mean component, $\xi_i$ are independent standard Gaussian random variables, $\lambda_i$ and $f_i(x)$ are eigenvalues and eigenfunctions of covariance functions described in Eq. (8), respectively. Here, the mean component is assumed to be 2. The 100 KL terms can preserve about 98% the field variance (defined as $\sum_{i=1}^{i=100} \lambda_i \Big/ \sum_{i=1}^{\infty} \lambda_i$).

Other settings for this case study, e.g., initial and boundary conditions, prior distributions of the location and release history, are the same with Case 1. Thus, there are 108 parameters for this case, i.e., 2 source location parameters ($x_s$, $y_s$), 6 source strength parameters ($S_{s1}, S_{s2}, S_{s3}, S_{s4}, S_{s5}, S_{s6}$), and 100 standard Gaussian variables for the truncated KL expansion ($\xi_1, \xi_2, \ldots, \xi_{100}$).

To infer these parameters, there are 40 sampling locations to provide hydraulic head at one time and concentration measurements every $2[\text{T}]$ from $t = 2[\text{T}]$ to $t = 10[\text{T}]$. The sampling locations are represented with blue dots in Figure 6(a). The concentration and hydraulic head measurements are generated with reference parameters with additive Gaussian errors $\varepsilon \sim N(0, 0.005^2)$ and $\varepsilon \sim N(0, 0.001^2)$, respectively. The reference conductivity is also generated by the truncated KL expansion (100 terms).

[Figure 6]

For the iIS algorithm, in each iteration, the number of parameter samples for updating was 400, and the 50 parameter samples with the smallest likelihood values were replaced with the inverse Gaussian process method as described in step (5), 2.2., which means that the total number of model evaluations was 450 in each iteration. Figure 7 shows the log Gaussian likelihood values of $\mathbf{En_d}$ (the ensemble of perturbed measurements) and $\mathbf{F}$ (the ensemble of model responses corresponding to the 400 parameter samples) at each iteration. At the first 6 iterations, the log likelihood values of $\mathbf{F}$ are far beyond the range of $\log Lik^{\mathbf{En_d}}$, which means that the distance between $\mathbf{F}$ and measurements $\mathbf{d}$ is large. At the 7th iteration, most of the log likelihood values of $\mathbf{F}$ are close to $\log Lik^{\mathbf{En_d}}$, and the iIS algorithm stops at the 9th iteration. In other words, the number of model evaluations for the iIS algorithm is $9 \times 450 = 4,050$. For the

DREAM$_{(ZS)}$ algorithm, 5 parallel chains and altogether 300,000 model evaluations were called, where Gaussian likelihood with heteroscedastic measurement errors was adopted. The trace plots of the 8 source parameters (location and release history) generated by the iIS algorithm and DREAM$_{(ZS)}$ are shown in Figure 8 and Figure 9, respectively. In Figure 8, for the iIS algorithm, as the 50 parameter samples replaced in each iteration are not shown, only 3,600 parameter samples are plotted in this figure. It is obvious that, the iIS algorithm can find the true source parameters much faster than MCMC.

[Figure 7]

[Figure 8]

[Figure 9]

For the conductivity field represented with 100 KL terms, Figure 7 (b) shows the true log K field. Using parameter samples with the biggest likelihood values for the iIS and the DREAM$_{(ZS)}$ algorithms, respectively, the log K field estimations are shown in Figure 7(c-d). It can be seen that, both algorithms can obtain log K fields close to the true reference log K field.

## 4. Conclusions

In this paper, an efficient approach entitled inverse iterative simulation is proposed to efficiently identify the contaminant source characteristics, together with the hydraulic conductivity field. The iIS algorithm utilizes a simple approach borrowed from Ensemble Smother (ES) to update model parameter and an inverse Gaussian process approach to improve the accuracy of parameter updating.

The efficiency and accuracy of the developed iIS algorithm in estimating contaminant source and conductivity field parameters were tested in two numerical case studies. In the first case study, with 8 contaminant source parameters and 3 hydraulic conductivity parameters, the iIS algorithm can obtain very close posterior parameter distributions compared with MCMC algorithm. In the second case study, with 8 contaminant source parameters and 100 KL expansion terms to represent the conductivity field, the iIS algorithm can obtain accurate estimation with very few model evaluations. Meanwhile, the time needed by the iIS algorithm could be further reduced through parallel computation.

## Acknowledgments.

Computer codes used are available upon request to the corresponding author.

We acknowledge Jasper Vrugt from University of California, Irvine for providing us with codes of DREAM$_{(ZS)}$.

## References

Evensen, G., and P. J. Van Leeuwen (2000), An ensemble Kalman smoother for nonlinear dynamics, *Mon. Weather Rev.*, *128*(6), 1852-1867, doi: 10.1175/1520-0493(2000)128<1852:AEKSFN>2.0.CO;2.

Haario, H., M. Laine, A. Mira, and E. Saksman (2006), DRAM: efficient adaptive MCMC, *Stat. Comput.*, *16*(4), 339-354, doi: 10.1007/s11222-006-9438-0.

Harbaugh, A. W. (2005), *MODFLOW-2005: The US Geological Survey Modular Ground-water Model--the Ground-water Flow Process*, U.S. Geol. Sur., Reston, VA, http://pubs.usgs.gov/tm/2005/tm6A16/.

Liu, C., and W. P. Ball (1999), Application of inverse methods to contaminant source identification from aquitard diffusion profiles at Dover AFB, Delaware, *Water Resour. Res.*, *35*(7), 1975-1985, doi: 10.1029/1999WR900092.

Mahinthakumar, G., and M. Sayeed (2005), Hybrid genetic algorithm—local search methods for solving groundwater source identification inverse problems, *J. Water Resour. Plan. Manage.*, *131*(1), 45-57, doi: 10.1061/(ASCE)0733-9496(2005)131:1(45).

Sidauruk, P., A. D. Cheng, and D. Ouazar (1998), Ground water contaminant source and transport parameter identification by correlation coefficient optimization,

*Groundwater*, *36*(2), 208-214, doi: 10.1111/j.1745-6584.1998.tb01085.x.

Snodgrass, M. F., and P. K. Kitanidis (1997), A geostatistical approach to contaminant source identification, *Water Resour. Res.*, *33*(4), 537-546, doi: 10.1029/96WR03753.

Sun, A. Y. (2007), A robust geostatistical approach to contaminant source identification, *Water Resour. Res.*, *43*(2), W02418, doi: 10.1029/2006WR005106.

Tartakovsky, D. M. (2013), Assessment and management of risk in subsurface hydrology: A review and perspective, *Adv. Water Resour.*, *51*, 247-260, doi: 10.1016/j.advwatres.2012.04.007.

Vrugt, J. A., C. J. Ter Braak, M. P. Clark, J. M. Hyman, and B. A. Robinson (2008), Treatment of input uncertainty in hydrologic modeling: Doing hydrology backward with Markov chain Monte Carlo simulation, *Water Resour. Res.*, *44*(12), W00B09, doi: 10.1029/2007wr006720.

Vrugt, J. A., C. Ter Braak, C. Diks, B. A. Robinson, J. M. Hyman, and D. Higdon (2009), Accelerating Markov chain Monte Carlo simulation by differential evolution with self-adaptive randomized subspace sampling, *Int. J. Nonlin. Sci. Num.*, *10*(3), 273-290, doi: 10.1515/IJNSNS.2009.10.3.273.

Wang, H., and X. Jin (2013), Characterization of groundwater contaminant source using Bayesian method, *Stoch. Env. Res. Risk A.*, *27*(4), 867-876, doi: 10.1007/s00477-012-0622-9.

Woodbury, A., E. Sudicky, T. J. Ulrych, and R. Ludwig (1998), Three-dimensional plume source reconstruction using minimum relative entropy inversion, *J. Contam. Hydrol.*, *32*(1), 131-158, doi: 10.1016/S0169-7722(97)00088-0.

Yeh, H. D., T. H. Chang, and Y. C. Lin (2007), Groundwater contaminant source identification by a hybrid heuristic approach, *Water Resour. Res.*, *43*(9), doi: 10.1029/2005WR004731.

Zeng, L., L. Shi, D. Zhang, and L. Wu (2012), A sparse grid based Bayesian method for contaminant source identification, *Adv. Water Resour.*, *37*, 1-9, doi: 10.1016/j.advwatres.2011.09.011.

Zhang, D., and Z. Lu (2004), An efficient, high-order perturbation approach for flow in random porous media via Karhunen–Loeve and polynomial expansions, *J. Comput.*

*Phys.*, *194*(2), 773-794, doi: 10.1016/j.jcp.2003.09.015.

Zhang, G., D. Lu, M. Ye, M. Gunzburger, and C. Webster (2013), An adaptive sparse-grid high-order stochastic collocation method for Bayesian inference in groundwater reactive transport modeling, *Water Resour. Res.*, *49*(10), 6871-6892, doi: 10.1002/wrcr.20467.

Zhang, J., L. Zeng, C. Chen, D. Chen, and L. Wu (2015), Efficient Bayesian experimental design for contaminant source identification, *Water Resour. Res.*, *51*(1), 576-598, doi: 10.1002/2014WR015740.

Zheng, C., and P. P. Wang (1999), MT3DMS: a modular three-dimensional multispecies transport model for simulation of advection, dispersion, and chemical reactions of contaminants in groundwater systems; documentation and user's guide*Rep.*, DTIC Document, http://www.geology.wisc.edu/courses/g727/mt3dmanual.pdf.

# Table captions:

**Table 1**. Prior range, true value, mean (Mean) and standard deviation (SD) values obtained by the iIS algorithm for each unknown parameter in Case Study 1.

# Figure captions:

**Figure 1.** Flow domain for Case Study 1.

**Figure 2.** Log likelihood values of model outputs ensemble and perturbed measurements ensemble in each iteration for Case Study 1.

**Figure 3.** Trace plots of (a, b) source location parameters, (c-h) source strength parameters and (i-k) log conductivity parameters in Case Study 1 generated by the iIS algorithm.

**Figure 4.** Probability distributions of contaminant transport model parameters inferred with DREAM$_{(ZS)}$ (represented by blue lines) and the iIS algorithm (represented by red lines). The true values are represented by vertical black lines.

**Figure 5.** Probability distributions of contaminant transport model parameters inferred with DREAM$_{(ZS)}$ (represented by blue lines), the iIS algorithm with iGP (represented by red lines) and the iIS algorithm without iGP (represented with purple lines). The true values are represented by vertical black lines.

**Figure 6.** (a) The flow domain and measurement locations (blue dots); (b) True log K field; (c) Log K field estimated with the iIS algorithm; (d) Log K field estimated with the DREAM$_{(ZS)}$ algorithm.

**Figure 7.** Log likelihood values of model outputs ensemble and perturbed measurements ensemble in each iteration for Case Study 2.

**Figure 8.** Trace plots of (a, b) source location parameters, (c-h) source strength parameters in Case Study 2 generated by the iIS algorithm.

**Figure 9.** Trace plots of (a, b) source location parameters, (c-h) source strength parameters in Case Study 2 generated by the DREAM$_{(ZS)}$ algorithm.

# Tables

Table 2. Prior range, true value, mean (Mean) and standard deviation (SD) values obtained by the iIS algorithm for each unknown parameter in Case Study 1.

|        | Range | True value | Mean  | SD      |
|--------|-------|------------|-------|---------|
| $x_s$    | [3 5] | 3.156      | 3.193 | 0.0210  |
| $y_s$    | [4 6] | 4.770      | 4.773 | 0.00190 |
| $S_{s1}$ | [0 8] | 6.239      | 5.933 | 0.118   |
| $S_{s2}$ | [0 8] | 5.667      | 5.899 | 0.0781  |
| $S_{s3}$ | [0 8] | 4.728      | 4.561 | 0.0789  |
| $S_{s4}$ | [0 8] | 3.016      | 3.256 | 0.0918  |
| $S_{s5}$ | [0 8] | 3.151      | 2.945 | 0.0727  |
| $S_{s6}$ | [0 8] | 3.427      | 3.621 | 0.0657  |
| $Y_1$    | [1 3] | 1.352      | 1.359 | 0.00890 |
| $Y_2$    | [1 3] | 2.722      | 2.758 | 0.0433  |
| $Y_3$    | [1 3] | 2.312      | 2.304 | 0.0174  |

# Figures

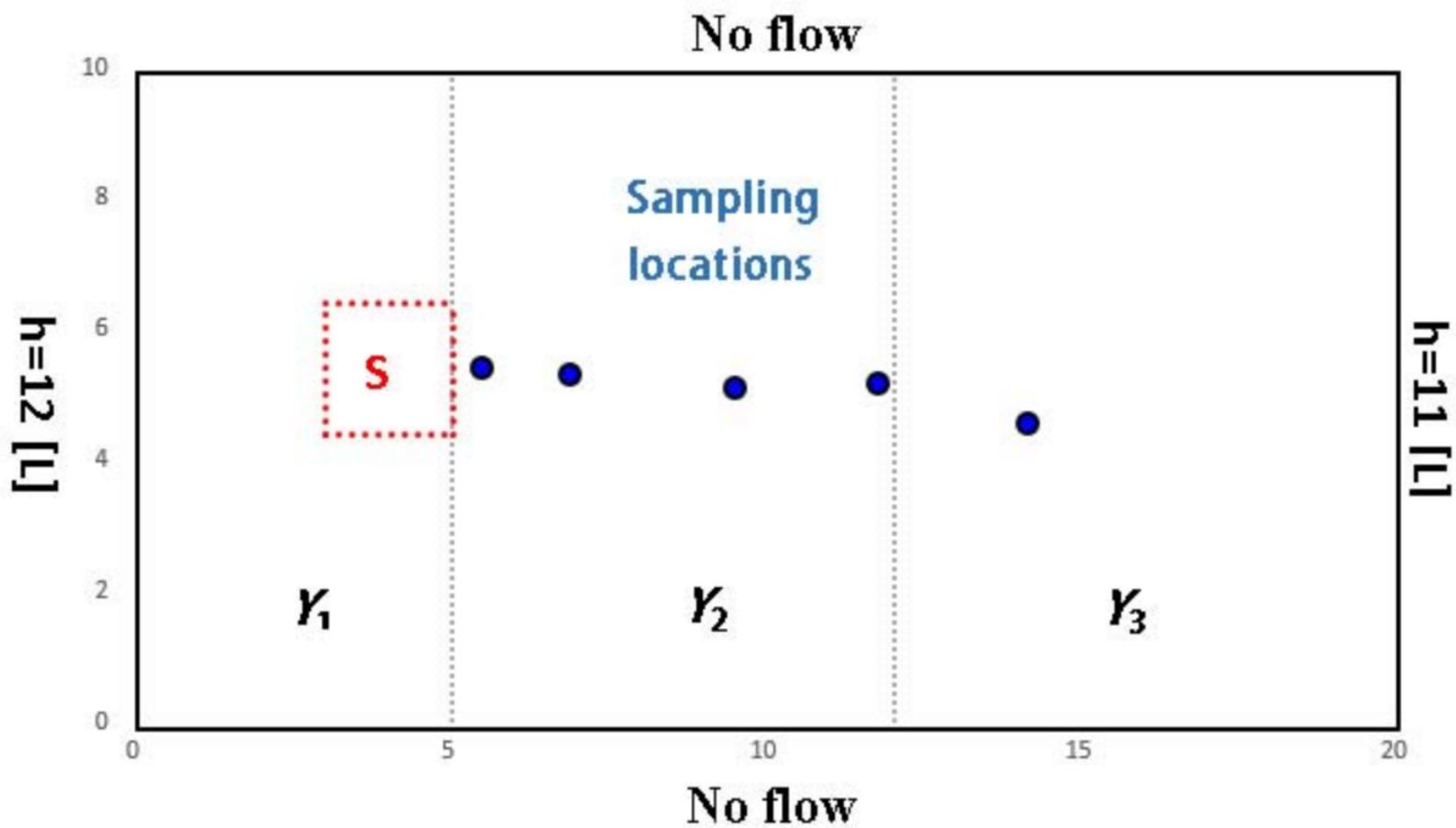

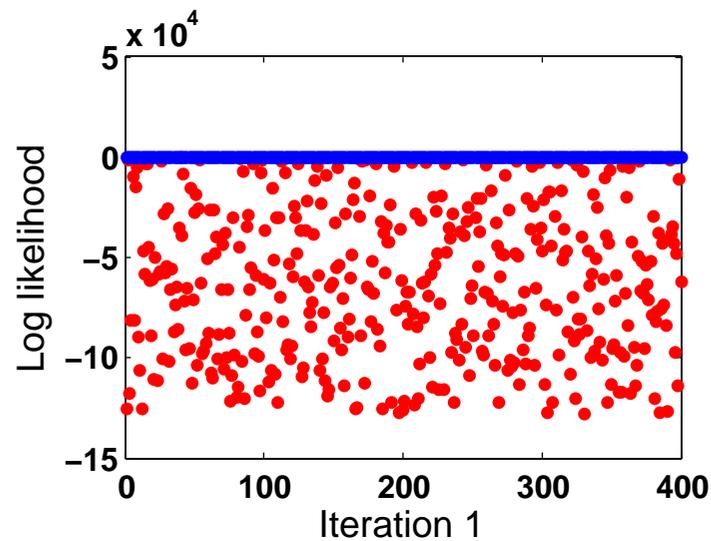 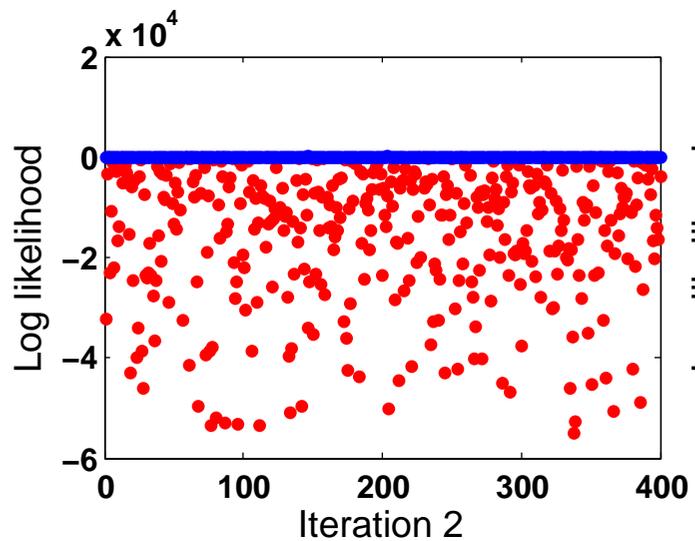 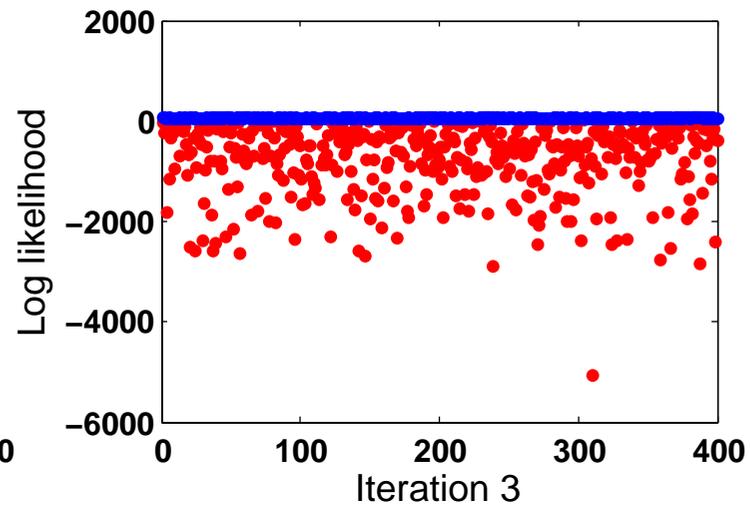
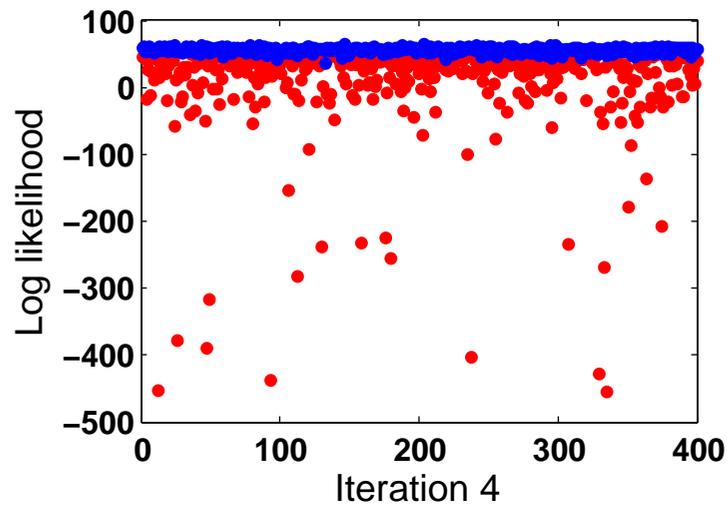 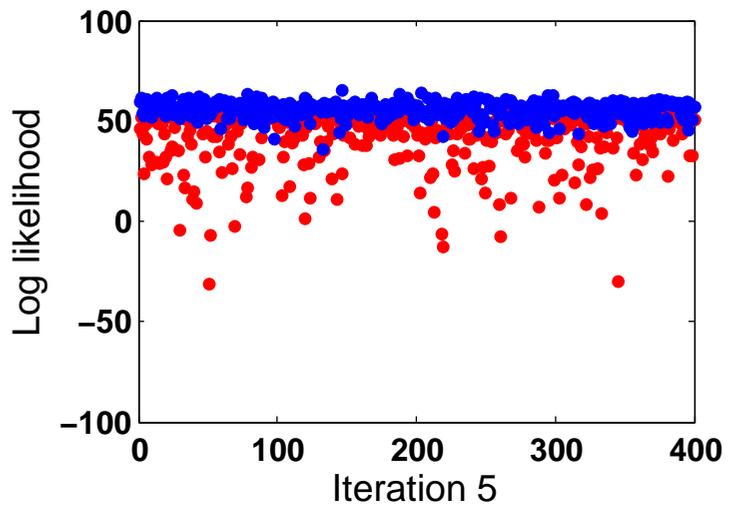 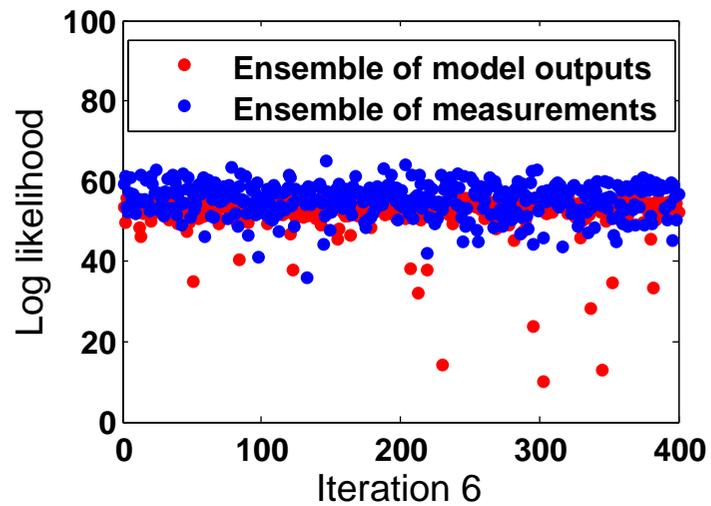

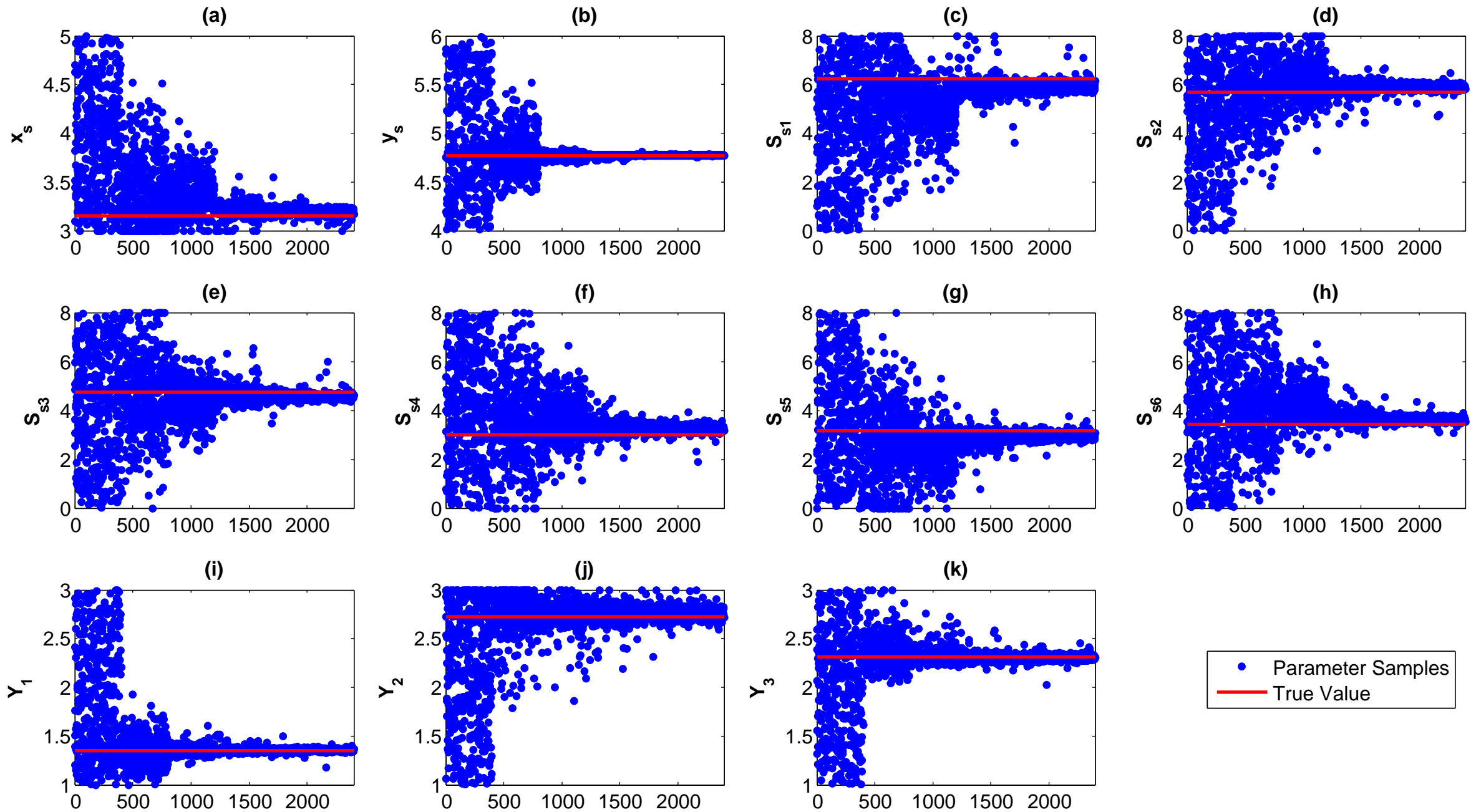

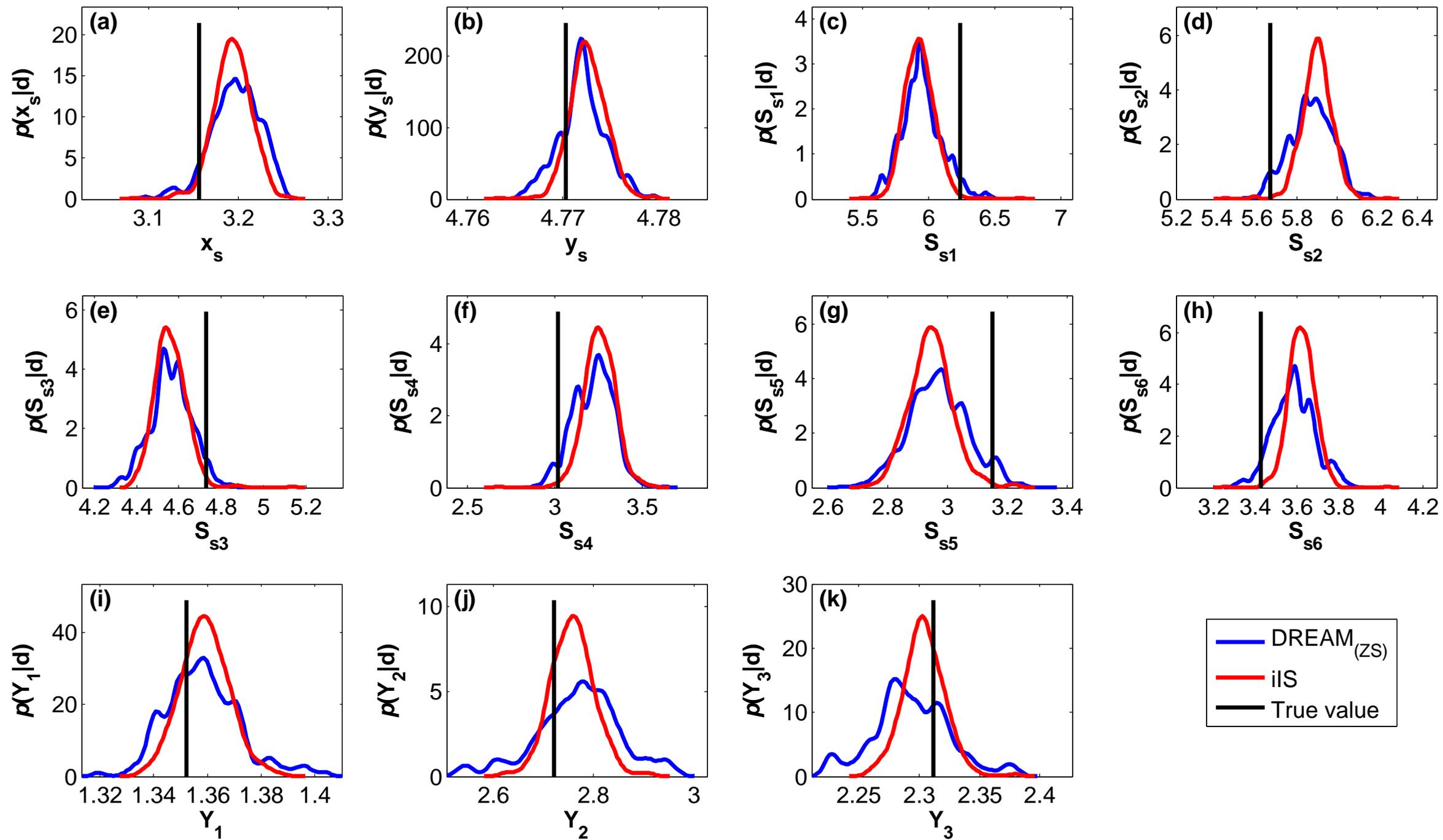

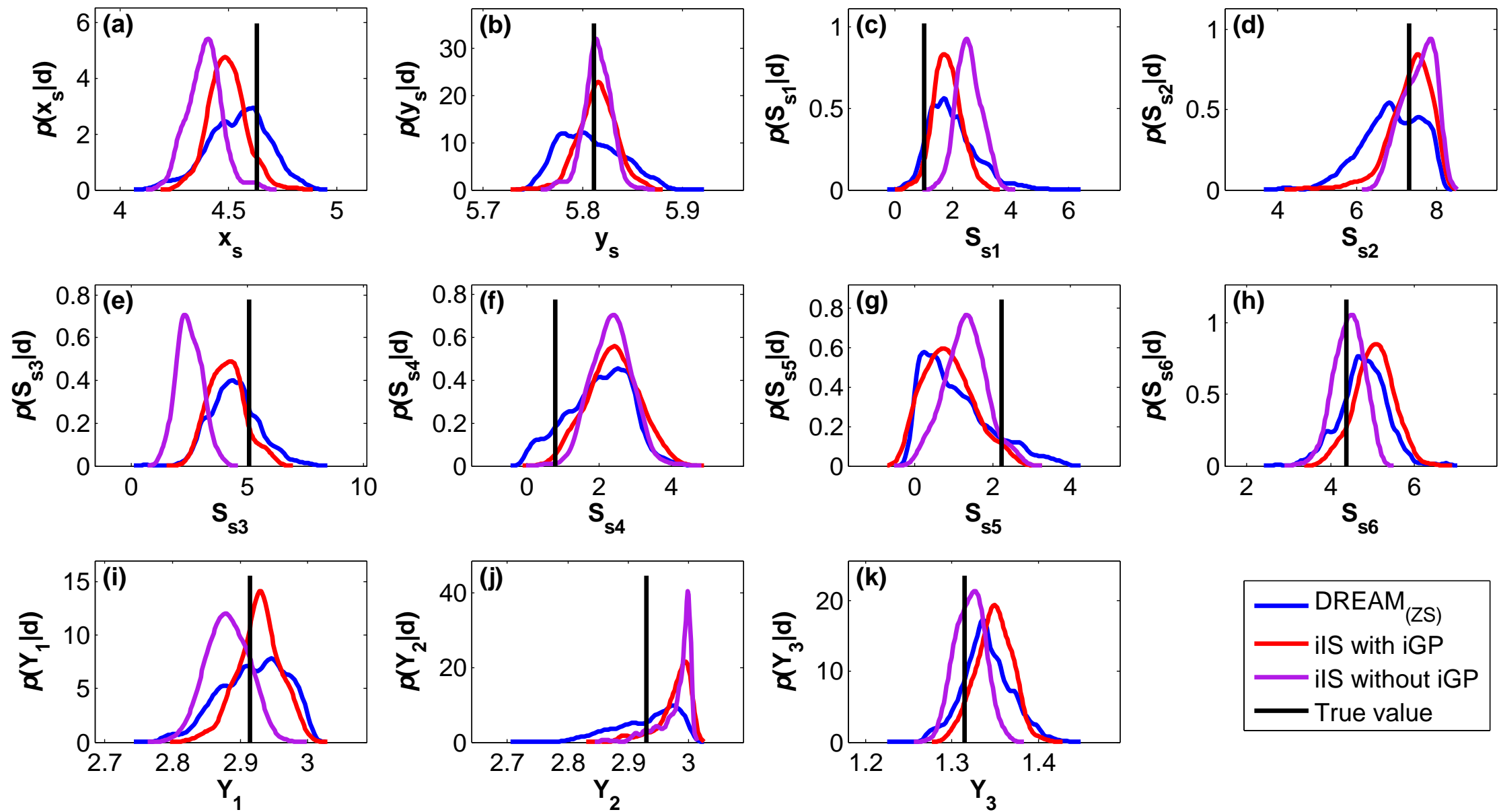

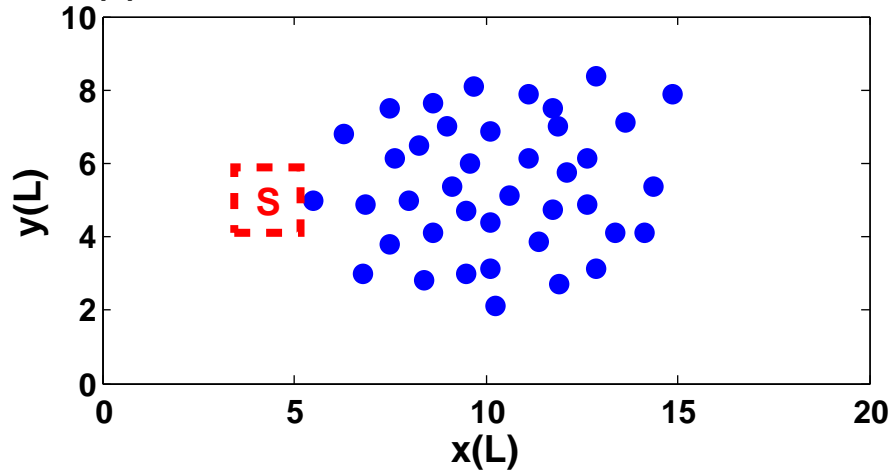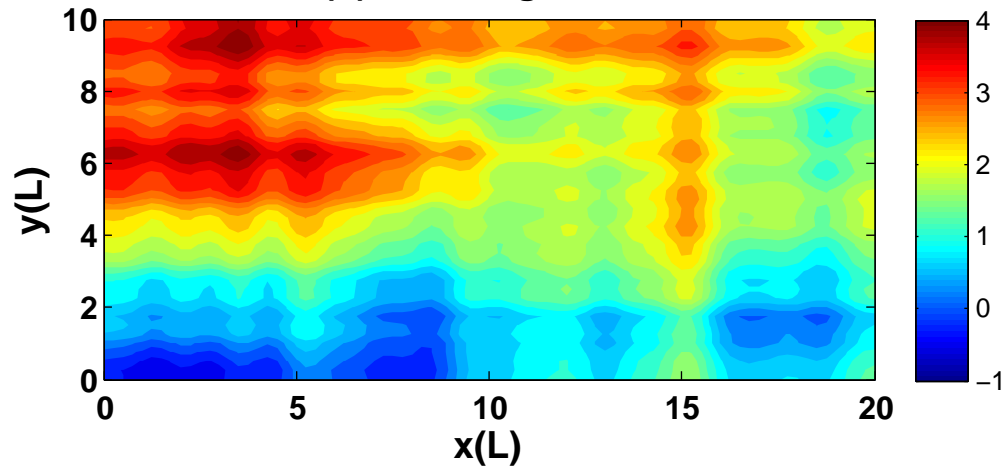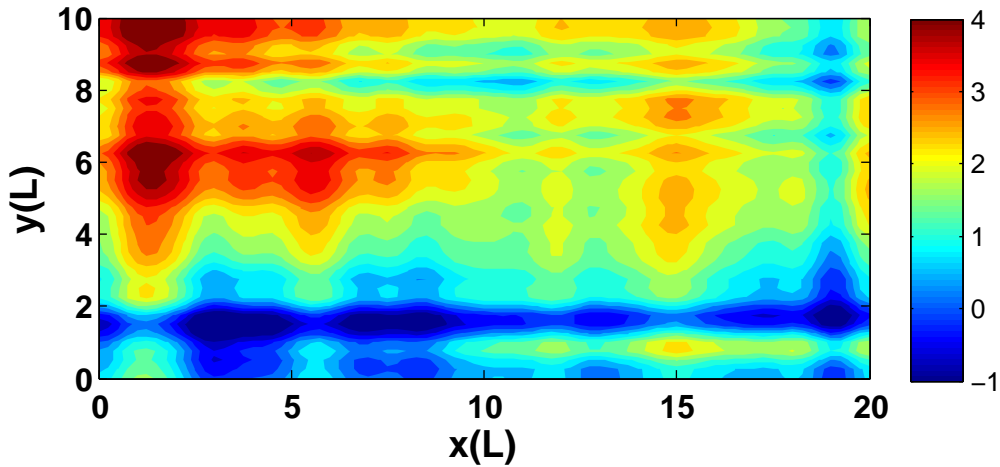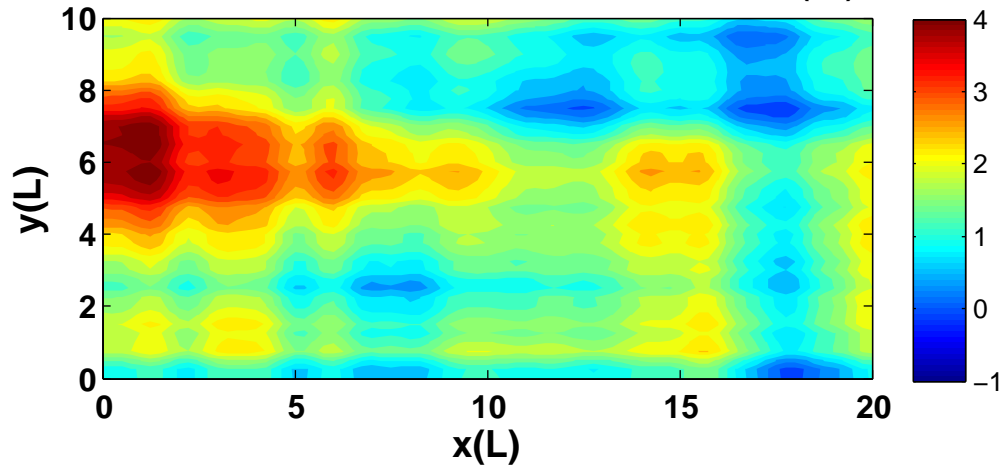

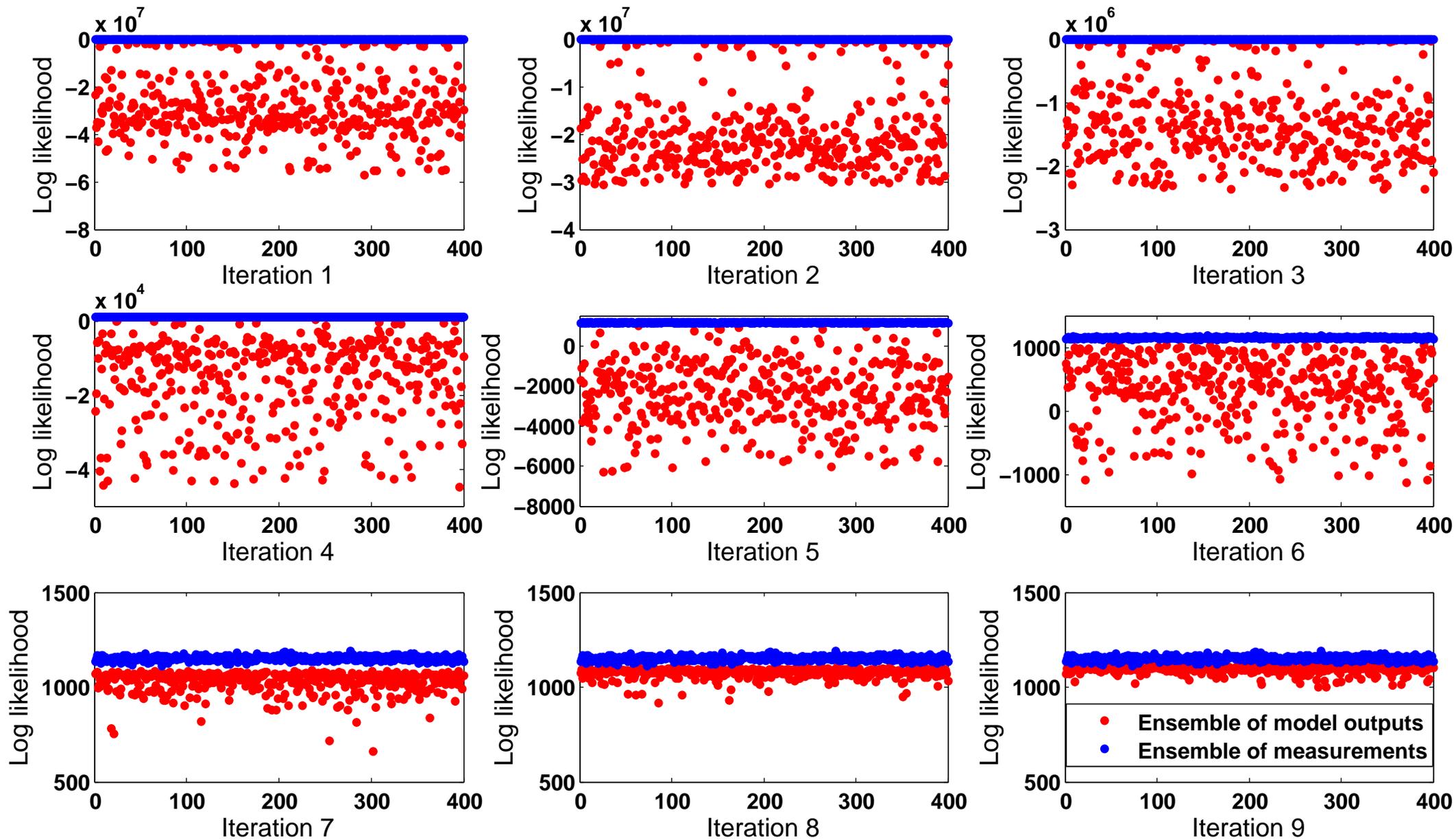

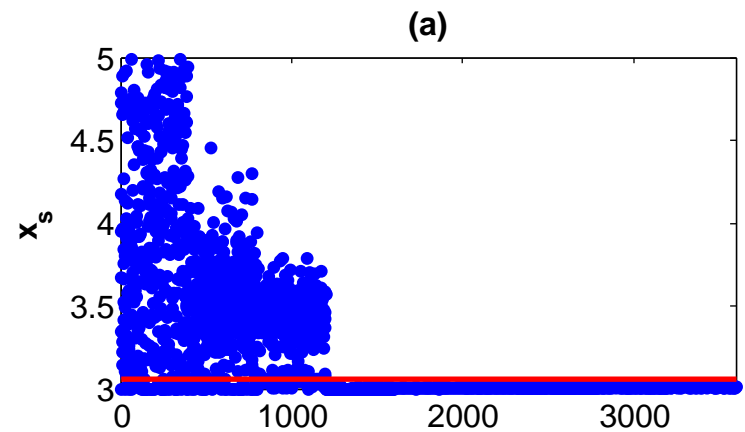
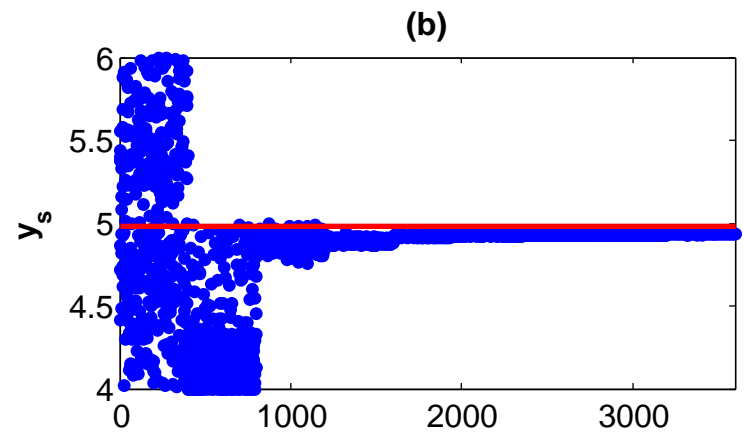
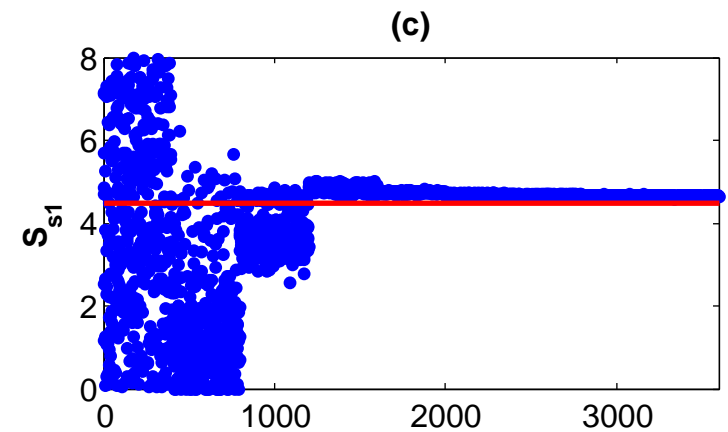
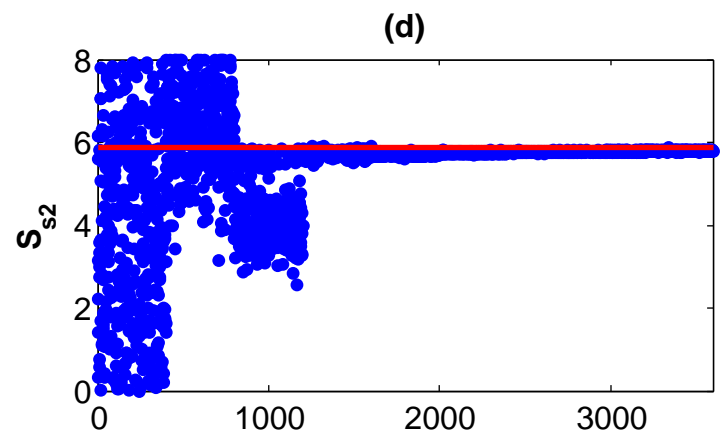
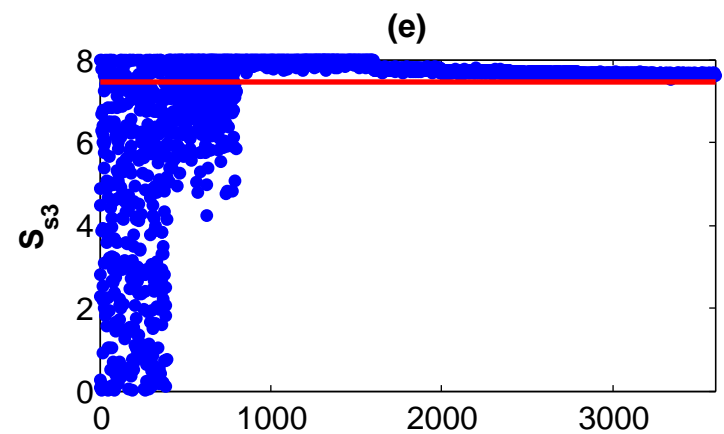
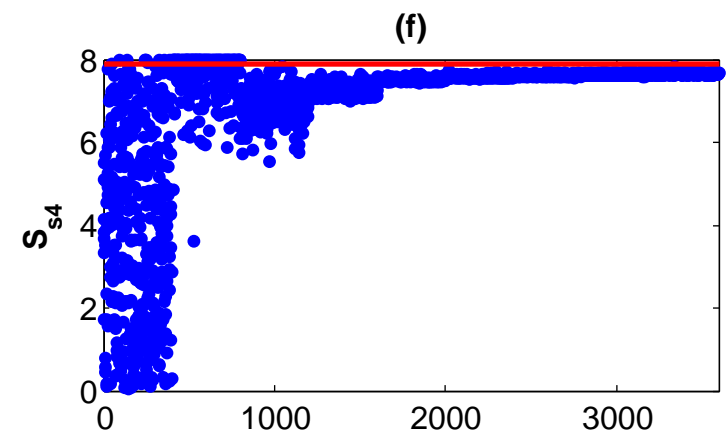
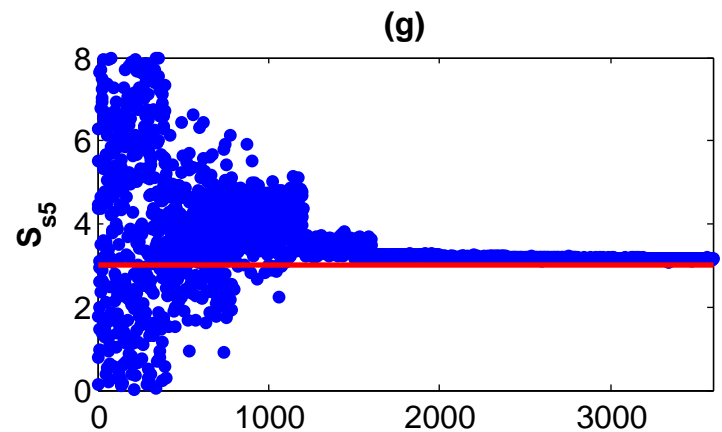
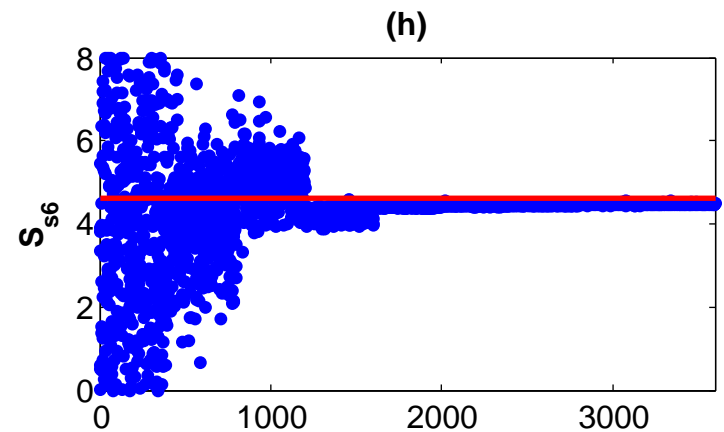
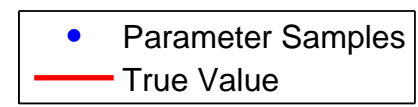

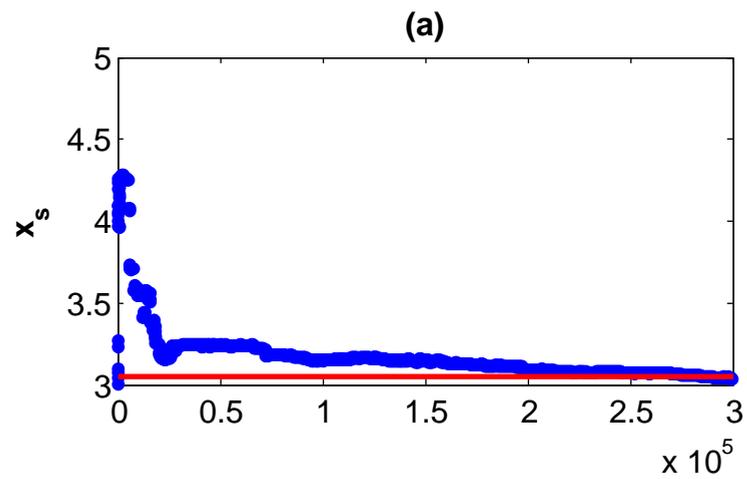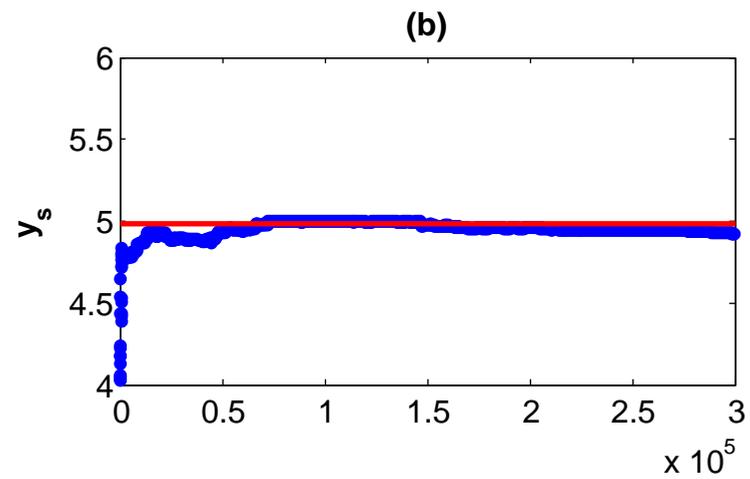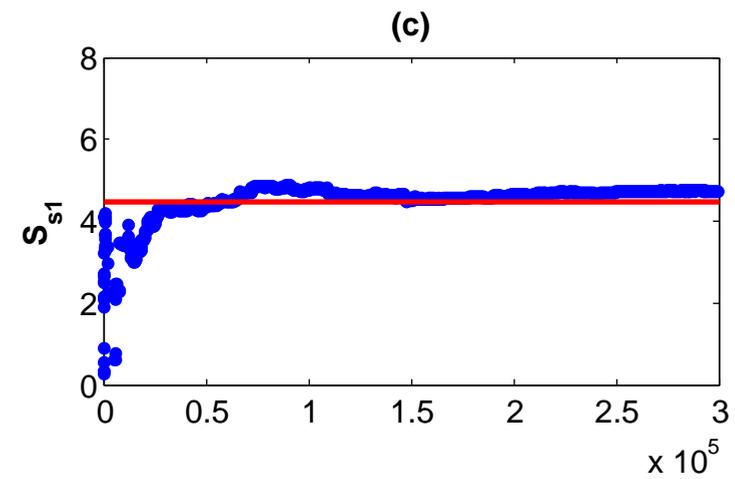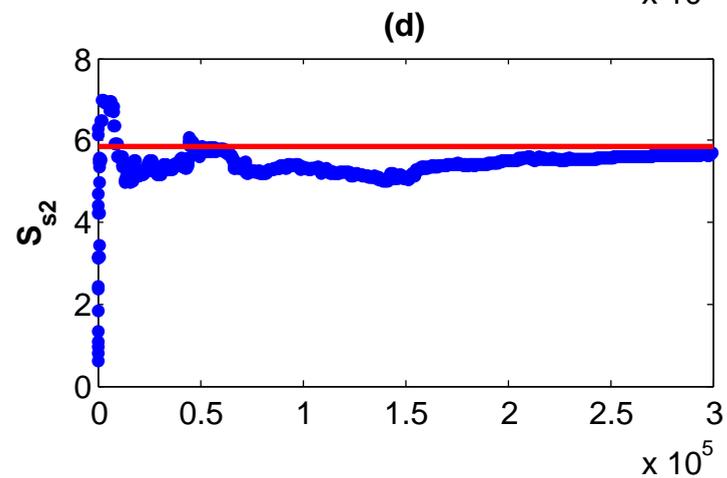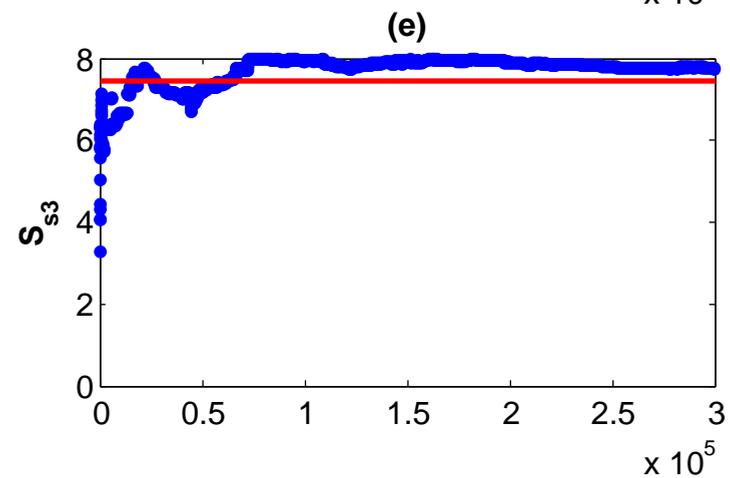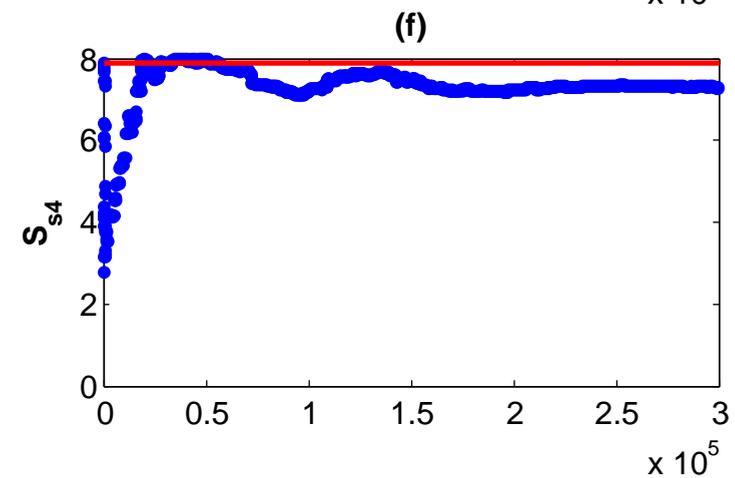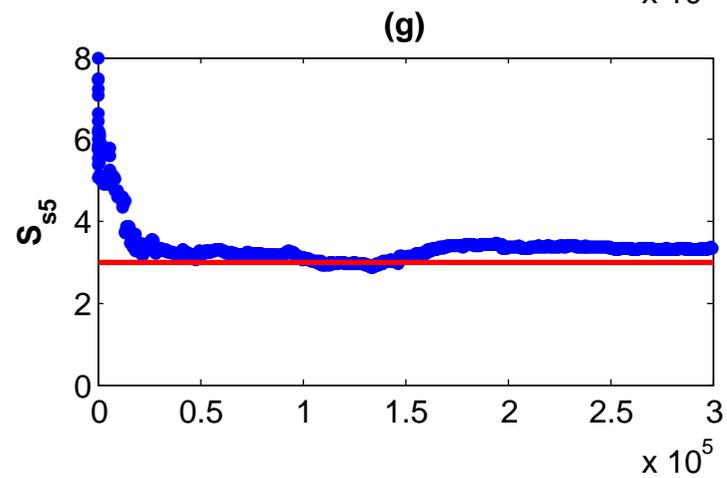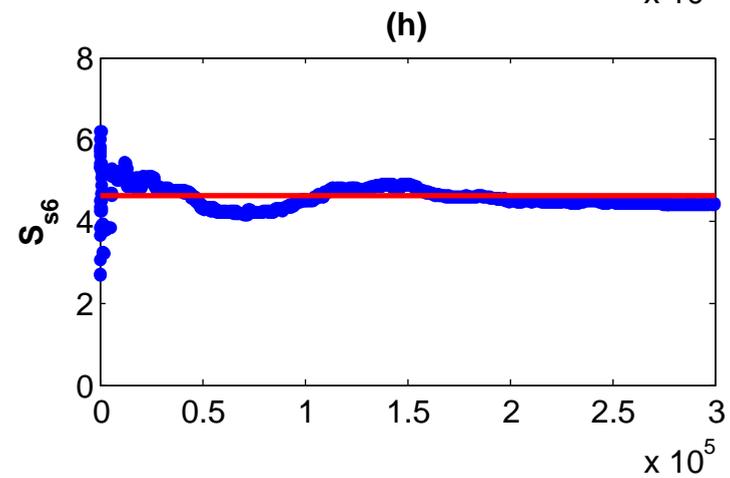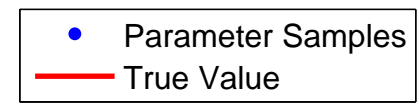